\documentclass{amsart}
\usepackage[T1]{fontenc}   
\usepackage{graphicx}
\usepackage{amsfonts,amsmath,amssymb,eucal,xypic}
\usepackage[all]{xy}
\usepackage{amsthm}

\begin{document}

\noindent                                             
\begin{picture}(150,36)                               
\put(5,20){\tiny{Submitted to}}                       
\put(5,7){\textbf{Topology Proceedings}}              
\put(0,0){\framebox(140,34){}}                        
\put(2,2){\framebox(136,30){}}                        
\end{picture}                                        
\vspace{0.5in}

\renewcommand{\bf}{\bfseries}
\renewcommand{\sc}{\scshape}
\vspace{0.5in}

\title[On the Diagonal Actions of $\mathbb{Z}$ on $\mathbb{Z}^n$]%
{On the Diagonal Actions of $\mathbb{Z}$ on $\mathbb{Z}^n$}

\author{Jordan Sahattchieve}
\address{Department of Mathematics;University of Michigan;
Ann Arbor, Michigan 48109}
\email{jantonov@umich.edu}
\thanks{The author was partially supported by NSF Research Training Grant DMS-0602191.}




\newtheorem{thm}{Theorem}[section]
\newtheorem{lem}[thm]{Lemma}
\newtheorem{pro}[thm]{Proposition}
\newtheorem{cor}[thm]{Corollary}
\newtheorem{ex}[thm]{Example}
\theoremstyle{defn}
\newtheorem{defn}[thm]{Definition}

\newtheorem{remark}[thm]{Remark}
\numberwithin{equation}{section}

\begin{abstract}
In this paper, we study the actions of $\mathbb{Z}$ on $\mathbb{Z}^n$ from a dynamical perspective.  The motivation for this study comes from the notion of bounded packing introduced by Hruska and Wise in \cite{HW}.  We also introduce the notion of coset growth for a finitely generated group.  Our analysis yields the following two results: bounded packing in certain semidirect products of $\mathbb{Z}^n$ with $\mathbb{Z}$ and a bound of the coset growth of the copy of $\mathbb{Z}$ on the right in $\mathbb{Z}^2 \rtimes \mathbb{Z}$ for the non-nilpotent groups of this type.
\end{abstract}

\maketitle

\section{\bf Introduction}
Our story begins with Sageev's celebrated cubing construction.  Recall that given a finitely generated group $G$ and a codimension-1 subgroup $H\leq G$, one constructs a CAT(0) cube complex $X_H$ on which $G$ acts \textit{essentially}, see \cite{Sag1}, \cite{Sag2} for details.  Naturally, this cubing depends on the choice of the codimension-1 subgroup $H$; however, it also depends on the choice of an \textit{invariant} $H$-\textit{set}.

In general, this cubing need not be finite dimensional, in which case the action of $G$ on $X_H$ has no hope of being cocompact.  The question of whether or not $X_H$ is finite dimensional leads us to the notion of \textit{bounded packing}.  Informally, $H$ has bounded packing in $G$, if given $R>0$, any sufficiently large collection of left cosets of $H$ in $G$ is guaranteed to contain a pair of cosets separated by a distance of at least $R$.

From the point of view of group actions on cube complexes, bounded packing is important because if $H\leq G$ has the property, it guarantees that for a suitably chosen invariant $H$-set $A$, Sageev's cubing $X_H$ which $A$ gives rise to will be finite dimensional.  This observation is implicit in \cite{Sag2} (see Lemma 3.2), even though the term \textit{bounded packing} itself does not appear there.  The notion was formally introduced by Hruska and Wise in \cite{HW} where they also establish the basic properties of bounded packing which we record in Section \ref{BPOverview}.

In this paper, we analyze the dynamics of the diagonal actions of $\mathbb{Z}$ on $\mathbb{Z}^n$.  Our interest in these actions came from the following question raised in \cite{HW}:\\\\
\textbf{Question} (Hruska, Wise \cite{HW}): Does every subgroup of a virtually polycyclic group have bounded packing?\\

The natural place to investigate this question from a geometric viewpoint is in the simplest non-abelian examples of polycyclic groups, namely the extensions of $\mathbb{Z}^n$ by $\mathbb{Z}$.  Our analysis of diagonal $\mathbb{Z}$-actions yields the following bounded packing result:\\\\
\textbf{Theorem 6.4.}  \textit{Let $\mathbb{Z}$ act on $\mathbb{Z}^n$ via a diagonalizable automorphism $\phi\in Aut(\mathbb{Z}^n)$ all of whose eigenvalues are real, then every subgroup of $\mathbb{Z}^n\rtimes_{\phi}\mathbb{Z}$ has bounded packing.}\\

Our techniques yield more than just this result on bounded packing - we are able to keep track of the growth of the set of left cosets of the copy of $\mathbb{Z}$ on the right in $\mathbb{Z}^2\rtimes\mathbb{Z}$.  We are thus naturally led to introduce \textit{coset growth} as a generalization of the notion of growth of groups and we use our estimates to obtain a bound for the coset growth of the copy of $\mathbb{Z}$ on the right in the non-nilpotent groups $\mathbb{Z}^2\rtimes\mathbb{Z}$:\\\\
\textbf{Theorem 7.5.}  \textit{Let $f_H(r)$ be the coset growth of $H=\left\langle t\right\rangle$ in $\mathbb{Z}^2\rtimes_{\phi}\mathbb{Z}$, where $\phi$ is a hyperbolic automorphism of $\mathbb{Z}^2$, then $f_H(r)\prec\alpha^r$.}\\

We would like to mention that the question on bounded packing in polycyclic groups was settled in its entirety by W Yang in \cite{YY}.  Our proof appeared in \cite{Sah1} shortly after Yang's preprint and uses completely different techniques.  Yang concludes that any separable subgroup has bounded packing in its ambient group and deduces bounded packing in polycyclic groups via a result of Malcev which establishes subgroup separability for these groups.  We did not appeal to Malcev's result on separability in polycyclic groups and were consequently able to obtain the estimate in Theorem \ref{Cogrowth}.

\section{\bf Bounded Packing - Definition and Basic Properties}\label{BPOverview}
If $G$ is a topological group, a metric $d$ on $G$ is called \textit{left-invariant} if for all $g\in G$, the map $a\mapsto ga$ is an isometry of $G$ to itself.  If $G$ is a finitely generated group, the word metric on $G$ corresponding to a finite generating set is clearly a left-invariant metric.  Throughout this section, $G$ will be a finitely generated group, $d$ will be the word metric unless otherwise stated, and $H$ will be a subgroup.
\begin{defn}\label{BP}(Bounded Packing) We say that $H$ has \textit{bounded packing} in $G$, if given any $r>0$ we can find a number $N$, which may depend on $G$, $H$, and $r$, such that any collection of $N$ distinct left cosets of $H$ in $G$ $\left\{g_1H,...,g_NH\right\}$ contains a pair $g_iH,g_jH$ which is separated by distance at least $r$, or $d(g_iH,g_jH)\geq r$.
\end{defn}

The original definition in \cite{HW} does not mention the word metric on $G$ specifically but rather refers to an arbitrary proper left-invariant metric.  It turns out that the choice of proper left-invariant metric on $G$ does not matter, since if $d_1$ and $d_2$ are two such metrics on $G$, the subgroup $H\subseteq G$ has bounded packing with respect to $d_1$ if and only if it has bounded packing with respect to $d_2$.  See Lemma 2.2 in \cite{HW}.

It follows immediately from Definition \ref{BP} that any finite index subgroup of $G$ has bounded packing in $G$ since there are no collections of $N$ distinct cosets for $N>|G:H|$.  Let us record this observation as:
\begin{lem}[Hruska, Wise \cite{HW}]\label{FinInd}Any finite index subgroup $K$ of a finitely generated countable group $G$ has bounded packing in $G$.
\end{lem}
Our next lemma concerns the behavior of the property of bounded packing under passing to subgroups or supergroups:
\begin{lem}[Hruska, Wise, \cite{HW}]\label{SubgrpPro}Suppose that $H\leq K\leq G$, and $G$ is a finitely generated group.
\begin{enumerate}
	\item If $H$ has bounded packing in $G$, then $H$ has bounded packing in $K$.
	\item If $H$ has bounded packing in $K$ and $K$ has bounded packing in $G$, then $H$ has bounded packing in $G$.
\end{enumerate}
\begin{proof}See Lemma 2.4 in \cite{HW}.
\end{proof}
\end{lem}
An important property of bounded packing is that it is preserved by replacing the original subgroup $H$ with a different subgroup $K\subseteq G$ which either contains it or is contained in it with finite index:
\begin{pro}[Hruska, Wise \cite{HW}]\label{Commensurability}Let $G$ be a finitely generated group.
\begin{enumerate}
	\item\label{Com1} If $H\leq K\leq G$ and $\left|K:H\right|<\infty$, then $H$ has bounded packing in $G$ if and only if $K$ has bounded packing in $G$.
	\item If $H\leq K\leq G$ and $\left|G:K\right|<\infty$, then $H$ has bounded packing in $K$ if and only if $H$ has bounded packing in $G$.
	\item If $H,K\leq G$ and $\left|G:K\right|<\infty$, then $H\cap K$ has bounded packing in $K$ if and only if $H$ has bounded packing in $G$.
\end{enumerate}
\begin{proof}See Proposition 2.5 in \cite{HW}.
\end{proof}
\end{pro}
Part (\ref{Com1}) of the proposition has the following immediate corollary:
\begin{cor}[Hruska, Wise \cite{HW}]Any finite subgroup of the finitely generated group $G$ has bounded packing.
\begin{proof}It suffices to show that the trivial group $\left\{1\right\}$ has bounded packing in $G$, but that is an easy consequence of the properness of the word metric on $G$.
\end{proof}
\end{cor}
The following two results summarize the elementary properties of bounded packing.  For us, Lemma \ref{ExactSeq} will be especially important.
\begin{lem}[Hruska, Wise, \cite{HW}]\label{ExactSeq}Let $1\longrightarrow N\longrightarrow G\longrightarrow Q\longrightarrow 1$ be a short exact sequence of groups where $G$ is finitely generated.  Let $H$ be a subgroup of $G$ which projects to the subgroup $\overline{H}$ of $Q$.  Then $\overline{H}$ has bounded packing in $Q$ if and only if $HN$ has bounded packing in $G$.
\begin{proof}See Lemma 2.8 in \cite{HW}.
\end{proof}
\end{lem}
\begin{cor}[Hruska, Wise \cite{HW}]\label{NormalSubgrp}Every normal subgroup $N$ of a finitely generated group $G$ has bounded packing.
\begin{proof}Follows immediately from \ref{ExactSeq} by taking $H=\left\{1\right\}$.
\end{proof}
\end{cor}
\begin{thm}[Hruska, Wise \cite{HW}]\label{Nilpotent}If $G$ is a finitely generated virtually nilpotent group, then each subgroup of $G$ has bounded packing in $G$.
\begin{proof}The proof is by induction on the length of the lower central series of $G$ and relies mostly on Lemma \ref{ExactSeq}.  See Theorem 2.12 in \cite{HW}.
\end{proof}
\end{thm}
The main result of \cite{HW} is about bounded packing in \textit{relatively hyperbolic} groups:
\begin{thm}[Hruska, Wise, \cite{HW}]\label{RelHypBP} Let $(G,\mathbb{P})$ be a relatively hyperbolic group, with a finite generating set $\mathcal{S}$, and let $H$ be a $\nu$-relatively quasiconvex subgroup of $G$.  Suppose that for each peripheral subgroup $P\in\mathbb{P}$ and each $g\in G$ the intersection $P\cap gHg^{-1}$ has bounded packing in $P$.  Then, $H$ has bounded packing in $G$.
\begin{proof}See Theorem 8.10 in \cite{HW}.
\end{proof}
\end{thm}
\section{\bf Polycyclic Groups and the Groups $P_{n,\phi}$}\label{PolyGrps}
First, let us define the main object of study:
\begin{defn}\label{PolyDef}A group $G$ is called \textit{polycyclic} if there exist subgroups $G_0=\left\{1\right\}\subseteq G_1\subseteq\cdots\subseteq G_{n-1}\subseteq G_n=G$ with $G_{i-1}$ normal in $G_i$, such that $G_i/G_{i-1}$ is a cyclic group.
\end{defn}
An easier to understand from geometrical perspective example of a polycyclic group is the semidirect product $\mathbb{Z}^n\rtimes_{\phi}\mathbb{Z}$, $\phi\in Aut(\mathbb{Z}^n)$.  As a set $\mathbb{Z}^n\rtimes_{\phi}\mathbb{Z}$ is just $\mathbb{Z}^{n+1}$, however $\mathbb{Z}^n\rtimes_{\phi}\mathbb{Z}$ is in general far from being quasi-isometric to $\mathbb{E}^{n+1}$.

Recall that the semidirect product $\mathbb{Z}^n\rtimes_{\phi}\mathbb{Z}$ has the following presentation:
\begin{center}
$\mathbb{Z}^n\rtimes_{\phi}\mathbb{Z}=\left\langle \mathbb{Z}^n,t:tat^{-1}=\phi(a),a\in\mathbb{Z}^n\right\rangle$,
\end{center}
where $t$ denotes a generator for the $\mathbb{Z}$ factor on the right.

As a matter of convenience, we shall denote the semidirect product of $\mathbb{Z}^n$ with $\mathbb{Z}$ by $P_{n,\phi}$.

The free abelian group $\mathbb{Z}^n$ is a normal subgroup of $P_{n,\phi}$ and the quotient $P_{n,\phi}/\mathbb{Z}^n$ is isomorphic to $\mathbb{Z}$.  Hence, we can write $P_{n,\phi}$ as a disjoint union of the left cosets of $\mathbb{Z}^n$ indexed by $\mathbb{Z}$: $P_{n,\phi}=\coprod_{i\in\mathbb{Z}}t^i\mathbb{Z}^n$.  We therefore conclude that we can construct the Cayley graph of $P_{n,\phi}$ in the following way: we take a union of $\mathbb{Z}$ copies of the Cayley graph of $\mathbb{Z}^n$ and insert an extra edge at each vertex joining $t^ia\in t^i\mathbb{Z}^n$ to $t^{i+1}\phi^{-1}(a)\in t^{i+1}\mathbb{Z}^n$, which corresponds to multiplication on the right by the generator $t$.

In Section \ref{BPinPn} below, we show that in a certain sense $P_{n,\phi}$ has only one interesting subgroup from the point of view of bounded packing, namely the cyclic subgroup $H=\left\langle t\right\rangle$.  Therefore, our first goal will be to show that $H$ has bounded packing in $P_{n,\phi}$.

It is clear from the description of the Cayley graph of $P_{n,\phi}$ that each left coset of $H$ is a transversal to the left cosets of $\mathbb{Z}^n$.  On the other hand, the obvious injection $\iota:\mathbb{Z}^n\rightarrow t^i\mathbb{Z}^n\subseteq P_{n,\phi}$ induces a graph monomorphism of the Cayley graphs $\iota:Cayley(\mathbb{Z}^n)\hookrightarrow Cayley(P_{n,\phi})$.  Therefore, in addition to the metric on $t^i\mathbb{Z}^n\subseteq P_{n,\phi}$ coming from the restriction of the word metric $d$ of $P_{n,\phi}$, we have a metric $d_{t^i\mathbb{Z}^n}$ which is the push-forward of the word metric on $\mathbb{Z}^n$.  Our first result shows that given any $D>0$, we can find a $N(D)\in\mathbb{N}$, such that any collection of $N(D)$ distinct left cosets of $H$ will contain a pair $aH,bH$ such that $d_{t^i\mathbb{Z}^n}(at^i,bt^i)>D$ for all $i\in\mathbb{Z}$, which in turn implies $d(aH,bH)>\dfrac{D}{2}$.  Then, having shown bounded packing for $H$, we prove bounded packing for all of the subgroups of $P_{n,\phi}$.

In order to show bounded packing of $H$ in $P_{n,\phi}$ as outlined above, we need to analyze the action of $\mathbb{Z}$ on $\mathbb{Z}^n$ via $\phi\in Aut(\mathbb{Z}^n)$.  This requires us to present some basic results on the dynamics of linear maps.

\section{\bf The Dynamics of Linear Maps}\label{DynLinMaps}
All of the material in this section is contained in the classical reference on dynamical systems by Hasselblatt and Katok \cite{KH}.  Their book, however, offers much more than we need, in fact, most of the material we require is contained in Section 1.2 of \cite{KH}.

Throughout this section, let $\phi:\mathbb{R}^n\rightarrow\mathbb{R}^n$ be a linear map.  It is well known that the algebraic behavior of a linear map $\phi$ is, to a large extent, determined by its \textit{eigenvalues}.  Just for the sake of completeness, we recall that $\lambda\in\mathbb{C}$ is an \textit{eigenvalue} for $\phi$, considered as a linear map $\mathbb{C}^n\rightarrow\mathbb{C}^n$, if there exists $v\in\mathbb{C}^n$ such that $\phi(v)=\lambda v$.  The set of eigenvalues of $\phi$ is called the spectrum of $\phi$, and the absolute value of the largest eigenvalue of $\phi$, denoted by $\rho(\phi)$, is called the \textit{spectral radius} of $\phi$.  In the case of bounded linear maps between Banach spaces, see the discussion below, the spectrum need not be finite.  It is, however, always a closed bounded and nonempty subset of $\mathbb{C}$, see for example \cite{La}.  For our purposes, we shall only need to consider finite spectra, in particular the spectra of linear maps on $\mathbb{R}^n$.  We denote the spectrum of a linear map $\phi$ by $Spec(\phi)$.

If we fix a basis for $\mathbb{R}^n$, such as the standard basis, we can represent the linear map $\phi$ by a matrix $A\in M_n(\mathbb{R})$.  To avoid repetition, let us fix a basis once and for all: In what follows, $\left\{e_i\right\}_{1\leq i\leq n}$ will denote the standard basis for $\mathbb{Z}^n$, $\mathbb{R}^n$, and $\mathbb{C}^n$, and the matrix representation of $\phi$ will always be assumed to be with respect to the standard basis unless otherwise stated.  We shall use $\phi$ and its matrix representation $A$ interchangeably, whenever there is little chance of confusion.
\begin{defn}The $n$-fold composition of $\phi$ with itself, $\phi\circ\phi\circ\ldots\circ\phi$ is called the $n$-th \textit{iterate} of $\phi$ and is denoted by $\phi^n$.
\end{defn}
We are interested in the sets $\mathcal{O}^+_v=\left\{\phi^i(v):i\in\mathbb{N}\right\}$, for $v\in\mathbb{R}^n$, and $\mathcal{O}_v=\left\{\phi^i(v):i\in\mathbb{Z}\right\}$, where we define $\phi^{-i}=(\phi^{-1})^n$ for $i>0$, whenever $\phi$ is invertible.  The former set is called the \textit{forward orbit} of $v$ under $\phi$, and the latter, simply the \textit{orbit} of $v$ under $\phi$.  In the case when $\phi$ is invertible, the forward orbit of $\phi^{-1}$ is called the \textit{backward orbit} of $\phi$.  It turns out that the eigenvalues of $\phi$ also determine the asymptotic behavior of the orbits of $\phi$, as the following proposition demonstrates:
\begin{pro}\label{ContrLinMap}For every $\delta>0$, there exists a norm on $\mathbb{R}^n$ such that\\ $\sup_{||v||=1}||\phi(v)||<\rho(\phi)+\delta$.
\begin{proof}See Proposition 1.2.2 in \cite{KH}.
\end{proof}
\end{pro}
It seems like Proposition \ref{ContrLinMap} is too vague to be useful, since it only establishes the bound on the norm of $\phi$ for \textbf{some} norm on $\mathbb{R}^n$.  However, it is a well know fact that all norms on a finite dimensional normed space are equivalent in the following strong sense:
\begin{pro}\label{EquivNorms}Let $\left\|\cdot\right\|_1$ and $\left\|\cdot\right\|_2$ be two norms on the finite dimensional vector space $V$, then there exists $C>0$ such that $\dfrac{1}{C}\left\|v\right\|_1\leq\left\|v\right\|_2\leq C\left\|v\right\|_1$, for all $v\in V$.
\begin{proof}Exercise.
\end{proof}
\end{pro}
Before we state the interesting for us corollary to Proposition \ref{ContrLinMap}, let us establish some notation.

The space of all linear maps on a normed vector space $V$ will be denoted by $\mathcal{L}(V)$.  It is easily seen that $\mathcal{L}(V)$ is itself a vector space under the usual rules for working with linear maps.  Given $\phi\in\mathcal{L}(V)$, the positive real number $\sup_{||v||=1}||\phi(v)||$ is called the \textit{norm} of the linear map $\phi$.  Unfortunately, however, it does not define a norm on all of $\mathcal{L}(V)$ since we can have a linear map on an infinite dimensional vector space $V$ for which this number is not finite.  We get around this by defining $\mathcal{B}(V)=\left\{\phi:V\rightarrow V:||\phi||<\infty\right\}$.  This subspace of $\mathcal{L}(V)$ is called the \textit{space of bounded linear maps on} $V$, and with the supremum norm defined above, $\mathcal{B}(V)$ becomes a Banach space.  A basic result of functional analysis asserts that a linear map is bounded if and only if it is continuous.  Actually, $\mathcal{B}(V)$ has more structure than just that of a Banach space, namely it possesses the structure of a \textit{Banach algebra} where multiplication is given by composition of linear maps.  It is an easy exercise to show that $|\phi|=\sup_{v\in V}\frac{||\phi(v)||}{||v||}$, in other words, $|\phi|\:||v||$ is an upper bound for $||\phi(v)||$ for \textbf{all} $v\in V$.  Even though we are interested only in the finite dimensional case, we have presented these basics on Banach spaces in order to put the present discussion into proper context.  We are now prepared to state the promised corollary to Proposition \ref{ContrLinMap}, which follows immediately in view of Proposition \ref{EquivNorms}:
\begin{cor}\label{ConToOrigin}If $\phi\in\mathcal{L}(\mathbb{R}^n)$ is a bounded linear map with spectral radius $\rho(\phi)<1$, then $\lim_{i\rightarrow\infty}\phi^i(v)=\textbf{0}$, for all $v\in\mathbb{R}^n$.  Equivalently, the only accumulation point of $\mathcal{O}^+(v)$ in $\mathbb{R}^n$ is $\textbf{0}$.  If in addition $\phi$ is also invertible, then the backward orbit of $\phi$ diverges to infinity.
\begin{proof}See Corollary 1.2.4 in \cite{KH}.
\end{proof}
\end{cor}
Even though we are not interested in the speed of convergence, the proof of the corollary shows that it is exponential.  Let us now consider the following easy to visualize example:
\begin{ex}Let $\phi:\mathbb{R}^2\rightarrow\mathbb{R}^2$ be the map given by the $2\times 2$ matrix
$
 \begin{pmatrix}
  \lambda & 0 \\
   0 & \lambda
 \end{pmatrix}
$, where $0<\lambda<1$.  The map $\phi$ act by multiplication by $\lambda<1$ and for every $v\in\mathbb{R}^2$, $\phi^i(v)$ converges to $\textbf{0}$ exponentially as $i\rightarrow\infty$.
\end{ex}
We are now ready to present our motivating example for much of the discussion that follows:
\begin{ex}\label{HypMtx}Let $0<\mu<1<\lambda$ and consider the map $\phi:\mathbb{R}^2\rightarrow\mathbb{R}^2$ given by
$A=
 \begin{pmatrix}
  \mu & 0 \\
   0 & \lambda
 \end{pmatrix}
$.  The matrix $A$ has an eigenbasis $\left\{e_1,e_2\right\}$ with eigenvalues $\mu$ and $\lambda$ respectively.  Let us define $E_-=\mathbb{R}e_1$ and $E_+=\mathbb{R}e_2$.  The forward orbit of any point in $E_-$ approaches the origin at an exponential speed, and the same is true for the backward orbit of any point in $E_+$.

Now, we proceed to analyze the behavior of the orbit of a point $(x_0,y_0)\in\mathbb{R}^2-(E_-\cup E_+)$ and to that end, without loss of generality, we assume that $x_0,y_0>0$.  We compute $\phi^i(x_0,y_0)=(\mu^ix_0,\lambda^iy_0)$, hence $\phi^i(x_0,y_0)$  lies on the curve $\mathcal{F}_{(x_0,y_0)}=\left\{(x,y)\in\mathbb{R}^2:xy^c=x_0y_0^c\right\}$, where the constant $c>0$ is defined by $\mu=\lambda^{-c}$.  Further, it is easily shown that $\phi$ leaves the curves $\mathcal{F}_{(x_0,y_0)}$ invariant.  In the special case where $\mu=\lambda^{-1}$, the curves $\mathcal{F}_{(x_0,y_0)}$ give a foliation of the plane by hyperbolae.
\end{ex}
Example \ref{HypMtx} motivates the following definition:
\begin{defn}(Hyperbolic map) A linear map $\phi:\mathbb{R}^n\rightarrow\mathbb{R}^n$ is called \textit{hyperbolic}, if the intersection of $Spec(\phi)$ with the unit circle in $\mathbb{C}$ is empty, or in other words, $\phi$ has no absolute value $1$ eigenvalues.
\end{defn}
Next, we recast some familiar notions from linear algebra from the point of view of dynamics.  Recall that to each $\lambda\in Spec(\phi)\cap\mathbb{R}$, for a given $\phi\in\mathcal{L}(\mathbb{R}^n)$, one associates its \textit{root space} $\bigcup_k Ker(\phi-\lambda I)^k$ denoted by $E_\lambda$.  Similarly, to each pair of complex conjugates $\lambda,\overline{\lambda}\in Spec(\phi)$, we define $E_{\lambda,\overline{\lambda}}$ to be the intersection of the root space of the complexified map $\phi\in\mathcal{L}(\mathbb{C}^n)$ with $\mathbb{R}^n$.  By analogy with Example \ref{HypMtx}, we set

\begin{equation}
E_-=\bigoplus_{|\lambda|<1}E_{\lambda}\oplus\bigoplus_{|\lambda|<1}E_{\lambda,\overline{\lambda}},
\end{equation} 

\begin{equation}
E_+=\bigoplus_{|\lambda|>1}E_{\lambda}\oplus\bigoplus_{|\lambda|>1}E_{\lambda,\overline{\lambda}},
\end{equation}

and

\begin{equation}
E_0=E_{-1}\oplus E_{1}\oplus\bigoplus_{|\lambda|=1} E_{\lambda,\overline{\lambda}}.
\end{equation}
Clearly then, $\mathbb{R}^n=E_-\oplus E_+\oplus E_0$, and also $\phi$ is a hyperbolic map if and only if $E_0=\left\{\textbf{0}\right\}$.

Everything we need from \cite{KH} is contained in the following:
\begin{pro}\label{RootSpDec}Let $\phi\in\mathcal{L}(\mathbb{R}^n)$ be a hyperbolic map.  Then:
\begin{enumerate}
	\item For every $v\in E_-$, the positive iterates $\phi^i(v)$ converge to the origin with exponential speed as $i\rightarrow\infty$, and if $\phi$ is invertible, then the negative iterates $\phi^i(v)$ go to infinity with exponential speed as $i\rightarrow -\infty$.
	\item For every $v\in E_+$, the positive iterates of $v$ go to infinity exponentially and if $\phi$ is invertible, then the negative iterates converge exponentially to the origin.
	\item For every $v\in\mathbb{R}^n-(E_-\cup E_+)$, the iterates $\phi^i(v)$ converge to infinity exponentially as $i\rightarrow\pm\infty$.
\end{enumerate}
\begin{proof}See Proposition 1.2.8 in \cite{KH}.
\end{proof}
\end{pro}

\section{\bf The Diagonal Actions of $\mathbb{Z}$ on $\mathbb{Z}^n$}\label{DiagAct}
This section contains a couple of technical lemmas which analyze the geometry of the space of left cosets of $H$ and lays a foundation for the work in Sections \ref{BPinPn} - \ref{CosetGr}.
 
Any automorphism of $\mathbb{Z}^n$ extends to an automorphism of $\mathbb{R}^n$.  On the other hand, an automorphism of $\mathbb{R}^n$ which preserves the integer lattice $\mathbb{Z}^n\subseteq\mathbb{R}^n$ restricts to an automorphism of $\mathbb{Z}^n$ if and only if it has determinant equal to $\pm 1$.  Thus, we identify $Aut(\mathbb{Z}^n)$ with the group of $n\times n$ integer matrices of determinant $\pm 1$.  This group is denoted by $GL_n(\mathbb{Z})$ and as we noted $GL_n(\mathbb{Z})\subseteq GL_n(\mathbb{R})$.

Now, given $\phi\in Aut(\mathbb{Z}^n)$ we turn our attention to the action of $\mathbb{Z}$ on $\mathbb{Z}^n$ induced by the homomorphism $\Phi:\mathbb{Z}\rightarrow Aut(\mathbb{Z}^n)$ given by $\Phi:1\mapsto\phi$.  If the automorphism $\phi$ is diagonalizable over $\mathbb{R}$, we shall say that $\mathbb{Z}$ \textit{acts diagonally} on $\mathbb{Z}^n$.  Our plan of attack on the problem of bounded packing of $H$ in $P_{n,\phi}$, as outlined at the closing of Section \ref{PolyGrps}, is to show that if one takes a large enough collection of distinct left cosets $a_1H,...,a_NH$, $a_i\in\mathbb{Z}^n$, then one is guaranteed that some pair will stay far apart in each left coset of $\mathbb{Z}^n$ in the metric $d_{t^i\mathbb{Z}^n}$.  It turns out that this is intimately connected with the geometry of the orbits of the $\mathbb{Z}$ action on $\mathbb{Z}^n$ via $\phi$.

We begin by showing that given a diagonal action of $\mathbb{Z}$ on $\mathbb{Z}^n$ by a hyperbolic automorphism, two translates of orbits of this action intersect in at most two points.  The reason for considering diagonal actions is simple, if all of the eigenvalues of $\phi$ are real and positive, the action of $\mathbb{Z}$ on $\mathbb{Z}^n$ extends to a \textit{flow} on $\mathbb{R}^n$ and we can treat the orbits of the $\mathbb{Z}$-action as differentiable curves in $\mathbb{R}^n$.
\begin{lem}\label{Orbits}Let $\mathbb{Z}$ act diagonally on $\mathbb{Z}^n$ via the hyperbolic map $\phi\in Aut(\mathbb{Z}^n)$ all of whose eigenvalues are positive.  If $\mathcal{O}_{z}$ and $\mathcal{O}_{w}$ are the orbits of the points $z,w\in\mathbb{Z}^n-\left\{0\right\}$, then $\left|\mathcal{O}_{z}\cap\left(a+\mathcal{O}_{w}\right)\right|\leq 2$ for any non-zero $a\in\mathbb{Z}^n$.
\begin{proof}
Let $\left\{v_1,...,v_n\right\}$ be an eigenbasis for $\phi$.  For $z\in\mathbb{Z}^n$, write $z=z_1v_1+...z_nv_n$, $a=a_1v_1+...+a_nv_n$.  Then, $\mathcal{O}_{z}$ is contained in the image of the curve $t\rightarrow (\lambda_1^tz_1,...,\lambda_n^tz_n)$.  Now, given $z,w\in\mathbb{Z}^n$, set $\alpha(t)=(\lambda_1^tz_1,...,\lambda_n^tz_n)$ and $\beta(t)=(\lambda_1^tw_1+a_1,...,\lambda_n^tw_n+a_n)$, and note that $\mathcal{O}_{z}\cap\left(a+\mathcal{O}_{w}\right)\subseteq Im(\alpha)\cap Im(\beta)$.  We show that there are at most two solutions $(t,t^{'})$ to the equation $\alpha(t)=\beta(t^{'})$.  Since $a\neq 0$, at least two coordinates of $a$ with respect to the chosen eigenbasis for $\phi$ are non-zero (otherwise $a$ would lie in the intersection of a 1-dimensional eigenspace of $\phi$ with $\mathbb{Z}^n$, which is impossible).  By reordering the basis, if necessary, we may assume that $a_1\neq 0$ and $a_2\neq 0$.  We consider the following cases:
\begin{enumerate}
	\item $z_1=0$ and $w_1\neq 0$:
	Considering the equation $\lambda_1^t z_1=\lambda_1^{t^{'}}w_1+a_1$, we see that there is at most one solution for $t^{'}$.  Now, we consider the subcases:
		\begin{itemize}
		\item $w_2=0$:
		Consider the equation $\lambda_2^t z_2=a_2$.  If $z_2=0$, this equation has no solution.  If $z_2\neq 0$, there exists at most one solution for $t$.
		\item $w_2\neq 0$ and $z_2\neq 0$:
		Considering $\lambda_2^t z_2=\lambda_2^{t^{'}} w_2+a_2$, with the value for $t^{'}$ we found above, we again conclude that there is at most one solution for $t$.
		\item $w_2\neq 0$ and $z_2=0$:
		In this case, since $z_1=z_2=0$, we must have some $z_i\neq 0$ since $z\neq 0$.  Assume that $z_3\neq 0$, and consider $\lambda_3^t z_3=\lambda_3^{t^{'}} w_3+a_3$.  With the value for $t^{'}$, there is again at most one solution for $t$.
		\end{itemize}
		\item $z_1\neq 0$, $w_1=0$, $z_2\neq 0$, $w_2\neq 0$:
		Consider $\lambda_1^t z_1=a_1$.  We conclude that there is at most one solution for $t$.  With this value for $t$ we solve $\lambda_2^t z_2=\lambda_2^{t^{'}} w_2+a_2$, and we find at most one solution for $t^{'}$.
		\item $z_1\neq 0$, $w_1\neq 0$, $z_2\neq 0$, $w_2\neq 0$:
		Consider the projections of the curves $\alpha$ and $\beta$ on the plane spanned by $\left\{v_1,v_2\right\}$.  These projections are the curves $t\rightarrow(\lambda_1^tz_1,\lambda_2^tz_2)$ and\\ $t\rightarrow(\lambda_1^tw_1+a_1,\lambda_2^tw_2+a_2)$ respectively.  Setting $t^{'}=t+k$, $\lambda_2=\lambda_1^p$, $\lambda_1^t=x$, $\lambda_1^k=y$, we obtain the system:\\\\ $\begin{array} {lcl} xz_1 & = & xyw_1+a_1 \\ x^pz_2 & = & x^py^pw_2+a_2 \end{array}$\\\\
Solving for $xy$ and substituting, we get:\\\\ $x^pz_2=\left(\frac{xz_1-a_1}{w_1}\right)^pw_2+a_2$.\\\\
We now show that this equation has at most 2 solutions for $x$.  Consider the function $f(x)=x^pz_2-\left(\frac{xz_1-a_1}{w_1}\right)^pw_2-a_2$.  Then,  $f^{'}(x)=p\left[z_2 x^{p-1}-\frac{z_1w_2}{w_1}\left(\frac{xz_1-a_1}{w_1}\right)^{p-1}\right]$.  Note that $p\neq 0$, as no eigenvalue of $\phi$ is equal to 1.  Note also that we may always arrange to have $p\neq 1$, since if $p=1$, then $\lambda_1=\lambda_2$, in which case we look for $a_i$, with $i\neq 1,2$, such that $a_i\neq 0$, and $\lambda_i\neq \lambda_1,\lambda_2$, and replace $v_2$ with $v_i$.  If no such $a_i$ exists, we conclude that $a$ lies in the intersection of $\mathbb{Z}^n$ and an eigenspace of $\phi$, which is impossible.  As $a_1\neq 0$, the only possible zero of $f^{'}$ must be located at $x=\frac{a_1}{z_1-w_1\left(\frac{z_2w_1}{z_1w_2}\right)^{\frac{1}{p-1}}}$.  Since $f$ has at most one more zero than $f^{'}$, we conclude that $f$ has at most two zeros.  Therefore, we have at most two solutions for $t$, and solving for $y$ concludes this final case.
\end{enumerate}
It is an easy exercise for the reader to convince himself that this exhausts all possible cases up to permutation of the variables involved.
\end{proof}
\end{lem}
Even though Lemma \ref{Orbits} was stated to handle the case where only one of the orbits was translated by $a\in\mathbb{Z}^n$, it takes no work at all to show that for any $z,w,a,b\in\mathbb{Z}^n$, we still have $|(a+\mathcal{O}_z)\cap(b+\mathcal{O}_w)|\leq 2$.  Now, using this counting argument, we prove the following interesting property of the diagonal action of $\mathbb{Z}$.  Consider the induced action of $\phi$ on the set $\mathcal{F}(\mathbb{Z}^n)$ of finite subsets of $\mathbb{Z}^n$.  We show that given $D>0$, there exist $N(D)\in\mathbb{N}$, such that any $\mathcal{C}\in\mathcal{F}(\mathbb{Z}^n)$ with $|\mathcal{C}|>N(D)$ is guaranteed to contain a pair of points whose images in every $\phi^k(\mathcal{C})$ are separated by a distance of at least $D$.  In other words, there exist $a,b\in\mathcal{C}$ such that $d(\phi^k(a),\phi^k(b))>D$ for all $k\in\mathbb{Z}$.
\begin{lem}\label{DistancesBddBelow}Let $\mathbb{Z}$ act diagonally on $\mathbb{Z}^n$ via the hyperbolic automorphism $\phi$ all of whose eigenvalues are positive.  Then given $D>0$, there exists $N>0$ such that for any collection of $m>N$ distinct points $a_1,a_2,...,a_m\in\mathbb{Z}^n$, there exist $i,j$ such that $d(\phi^k(a_i),\phi^k(a_j))>D$ for all $k\in\mathbb{Z}$.
\begin{proof}Let $S$ be the set $\overline{B}(0,D)\cap\mathbb{Z}^n$.  If $\left\|\phi^k(z_1)-\phi^k(z_2)\right\|\leq D$ for some $z_1,z_2\in\mathbb{Z}^n$ and $k\in\mathbb{Z}$, then $z_1-z_2$ is in the orbit $\mathcal{O}_z$ of some $z\in S$ under $\phi$.  Hence, $z_2\in z_1+\bigcup_{z\in S}\mathcal{O}_z$.  Now, if we have $\left\|\phi^k(a_i)-\phi^k(a_j)\right\|\leq D$ for all $1\leq i,j\leq m$ , then $a_m\in \bigcap_{1\leq l\leq m-1} (a_l+\bigcup_{z\in S}\mathcal{O}_z)$.  To finish the proof, we observe that Lemma \ref{Orbits} shows that $\left|\bigcap_{1\leq l\leq m-1} (a_l+\bigcup_{z\in S}\mathcal{O}_z)\right|\leq 2\left|S\right|^2\sim 2\left(Vol(\overline{B}(0,D))\right)^2$.
\end{proof}
\end{lem}
\section{\bf Bounded Packing in $P_{n,\phi}$}\label{BPinPn}
We now use the dynamical properties of the diagonal $\mathbb{Z}$-action on $\mathcal{F}(\mathbb{Z}^n)$ to prove bounded packing in $\mathbb{Z}^n\rtimes\mathbb{Z}$.
\begin{lem}\label{BddPackoft}Let $\mathbb{Z}$ act diagonally on $\mathbb{Z}^n$ via the hyperbolic automorphism $\phi$, all of whose eigenvalues are positive.  Let $P_{n,\phi}=\mathbb{Z}^n\rtimes_{\phi}\mathbb{Z}$ and let $t$ denotes a generator for the copy of $\mathbb{Z}$ in $G$ on the right, then the subgroup $H=\left<t\right>$ has bounded packing in $P_{n,\phi}$.
\begin{proof}Let $d$ denote the word metric on $G$, and let $d_{t^i\mathbb{Z}^n}$ denote the inherited Euclidean metric on $t^i\mathbb{Z}^n$ described in Section \ref{PolyGrps}.  We show that given $R>0$, there exists $D>0$ such that if $x,y\in t^i\mathbb{Z}^n$ and $d_{t^i\mathbb{Z}^n}(x,y)>D$, then $d(x,y)>R$.  After translating, we need only show that if $d_{\mathbb{Z}^n}(e,y)>D$, then $d(e,y)>R$.  Supposing to the contrary, that given a fixed $R>0$, for any $D>0$ there exists $y\in\mathbb{Z}^n$ with $d_{\mathbb{Z}^n}(e,y)>D$ and $d(e,y)<R$, we construct a sequence $\left\{y_k\right\}$ with $d_{\mathbb{Z}^n}(e,y_k)>k$ and $d(e,y_k)<R$.  Because the metric $d$ is proper, this is a contradiction.  Now, since $\mathbb{R}^n$ is quasi-isometric to the lattice $\mathbb{Z}^n$, Lemma \ref{DistancesBddBelow} shows that given $D>0$, there exists $N\in\mathbb{N}$ such that any collection of $m$ distinct cosets of $H$, with $m>N$, contains a pair, say $aH$, $bH$, such that in any coset $t^i\mathbb{Z}^n$, we have $d_{t^i\mathbb{Z}^n}(at^i, bt^i)>D$ for all $i\in\mathbb{Z}$.  This along with the discussion above shows that given $D>0$ we can find $N'\in\mathbb{N}$ such that any collection of $m$ distinct cosets, with $m>N'$, contains a pair, $aH$, $bH$, such that $d(at^i, bt^i)>D$ for all $i\in\mathbb{Z}$.  Finally, we show that this implies that $d(aH,bH)\geq \dfrac{D}{2}$.  To this end, let $at^i\in aH$, and $bt^j\in bH$.  We have two cases: either $|i-j|\geq \dfrac{D}{2}$, or $|i-j|<\dfrac{D}{2}$.  In the first case, it is clear that $d(at^i,bt^j)\geq |i-j|\geq \dfrac{D}{2}$.  In the second, we have: $d(at^i,bt^j)\geq|d(at^i,at^j)-d(at^j,bt^j)|>\left|\dfrac{D}{2}-D\right|=\dfrac{D}{2}$.
\end{proof}
\end{lem}
Having established bounded packing of $H$ in $P_{n,\phi}$ where the diagonal $\mathbb{Z}$-action on $\mathbb{Z}^n$ is via a hyperbolic automorphism, we now work to relax the hyperbolicity assumption.
\begin{lem}\label{RelaxHyp}Let $\mathbb{Z}$ act diagonally on $\mathbb{Z}^n$ via the automorphism $\phi$ all of whose eigenvalues are positive.  If $P_{n,\phi}=\mathbb{Z}^n\rtimes_{\phi}\mathbb{Z}$ and if $t$ denotes a generator for the right copy of $\mathbb{Z}$ in $P_{n,\phi}$, then the subgroup $H=\left<t\right>$ has bounded packing in $P_{n,\phi}$.
\begin{proof}The automorphism $\phi$ extends to an automorphism $\phi:\mathbb{R}^n\mapsto\mathbb{R}^n$ which has determinant equal to $\pm 1$.  We have $\mathbb{R}^n=\mathcal{H}\oplus E_1$, where $\mathcal{H}=E_-\oplus E_+$ is the sum of the expanding and contracting subspaces of $\phi$, and $E_1$ is the eigenspace for the eigenvalue $\lambda=1$ of $\phi$.  Note that since both $\mathcal{H}$ and $E_1$ are $\phi$-stable, so is $Z_{\mathcal{H}}=\mathbb{Z}^n\cap\mathcal{H}$, so $Z_{\mathcal{H}}\triangleleft P_{n,\phi}$.  Because $Z_{\mathcal{H}}=\mathbb{Z}^n\cap\mathcal{H}$ is a free abelian group on which $\phi$ acts as a hyperbolic automorphism, $\left<t\right>$ has bounded packing in $Z_{\mathcal{H}}\rtimes_{\phi}\mathbb{Z}=Z_{\mathcal{H}}\left<t\right>\subseteq P_{n,\phi}$ by Lemma \ref{BddPackoft}.  Now, $P_{n,\phi}/Z_{\mathcal{H}}$ embeds in $\left(\mathbb{R}^n\rtimes_{\phi}\mathbb{R}\right)/\mathcal{H}\cong\left(\mathbb{R}^n/\mathcal{H}\right)\rtimes_{\overline{\phi}}\mathbb{R}$
which is isomorphic to $\mathbb{R}^k$ since $\overline{\phi}$ is the identity.  Therefore, $P_{n,\phi}/Z_{\mathcal{H}}$ is abelian hence it has bounded packing with respect to all of its subgroups.  In particular, the image of $\left<t\right>$ has bounded packing in $P_{n,\phi}/Z_{\mathcal{H}}$, and therefore $Z_{\mathcal{H}}\left<t\right>$ has bounded packing in $P_{n,\phi}$.  As $\left<t\right>\subseteq Z_{\mathcal{H}}\left<t\right>\subseteq P_{n,\phi}$, we conclude that $\left<t\right>$ has bounded packing in $P_{n,\phi}$.
\end{proof}
\end{lem}
An immediate consequence of Lemma \ref{RelaxHyp} is:
\begin{cor}\label{BddPackofzt}With the notation from Lemma \ref{RelaxHyp}, for every $z\in\mathbb{Z}^n$, the subgroup $\left\langle zx\right\rangle$ has bounded packing in $P_{n,\phi}$.
\begin{proof}Since for any $z\in\mathbb{Z}^n$, we have:\\ $zxw\left(zx\right)^{-1}=z\phi(w)z^{-1}=\phi(w)=xzx^{-1}$, we have an automorphism of $P_{n,\phi}$ which is the identity on $\mathbb{Z}^n$ and which sends $x$ to $zx$.  The conclusion now follows immediately from Lemma \ref{RelaxHyp}.
\end{proof}
\end{cor}
This corollary shows that any cyclic subgroup of $P_{n,\phi}$ has bounded packing.

We finally have the means to prove the main result of this section:
\begin{thm}\label{BoundedPacking}Let $\mathbb{Z}$ act diagonally on $\mathbb{Z}^n$ via $\phi\in Aut(\mathbb{Z}^n)$ all of whose eigenvalues are real, then every subgroup of $P_{n,\phi}=\mathbb{Z}^n\rtimes_{\phi}\mathbb{Z}$ has bounded packing in $P_{n,\phi}$.
\begin{proof}By passing to a subgroup of finite index, we may assume that all of the eigenvalues of $\phi$ are real and positive.  Let $H\subseteq P_{n,\phi}$.  By further passing to a finite index subgroup of $P_{n,\phi}$, we may assume that the projection of $H$ to $\mathbb{Z}$ is onto, so that we have the exact sequence $1\longrightarrow W\longrightarrow H\longrightarrow\mathbb{Z}\longrightarrow 1$, where $W=\mathbb{Z}^n\cap H$.  Note that $W$ is normal in $P_{n,\phi}$: let $zt\in H$ be an element of $H$ which projects to a generator of the right-side copy of $\mathbb{Z}$ in $P_{n,\phi}$, where $z\in\mathbb{Z}^n$, then $ztWt^{-1}z^{-1}=W$, so $tWt^{-1}=z^{-1}Wz=W$.  Next, by passing to a subgroup of $P_{n,\phi}$ containing $H$ with finite index, we may assume that $\mathbb{Z}^n/W$ is free abelian.  Here is how: let $\mathbb{Z}^n/W=F\oplus T$, where $F$ is free abelian, and $T$ is torsion; the preimage of $T$ in $\mathbb{Z}^n$, say $T'$, under the quotient map is $\phi$-stable, $\left|T':W\right|<\infty$, so $T'\rtimes\mathbb{Z}$ is the desired subgroup.  Now, consider $G=P_{n,\phi}/W$.  Because the diagram\\\\ 
    $\xymatrix{ 1 \ar[r] & \mathbb{Z}^n \ar[r] \ar[d] & P_{n,\phi} \ar[r] \ar[d] & \mathbb{Z} \ar[r] & 1 \\
                1 \ar[r] & \mathbb{Z}^n/W \ar[r] & P_{n,\phi}/W \ar[r] & \mathbb{Z} \ar[r] & 1}$
\\\\\\
commutes, $G\cong(\mathbb{Z}^n/W)\rtimes_{\overline{\phi}}\mathbb{Z}\cong\mathbb{Z}^m\rtimes_{\overline{\phi}}\mathbb{Z}$, where $\overline{\phi}$ is the induced automorphism on the quotient $\mathbb{Z}^n/W$.  Corollary \ref{BddPackofzt} shows that the image of $\left<zt\right>$ has bounded packing in $P_{n,\phi}/W$, hence $H=W\left<zt\right>$ has bounded packing in $P_{n,\phi}$ by Lemma \ref{ExactSeq}.
\end{proof}
\end{thm}
Note that Theorem \ref{BoundedPacking} solves the question of bounded packing in polycyclic groups of length $\leq 3$ as the following corollary demonstrates:
\begin{cor}\label{Z2byZ}If $G$ is a polycyclic group of length $\leq 3$, $G$ has bounded packing with respect to all of its subgroups.
\begin{proof}The statement is trivially true if the length of $G$ is equal to 1 or 2, since in both cases $G$ is virtually abelian.  If $G$ has length 3, then $G$ contains a subgroup isomorphic to $\mathbb{Z}^2\rtimes_{\phi}\mathbb{Z}$ of finite index.  Unless $\phi$ is hyperbolic, $G$ will be virtually nilpotent; bounded packing in this case follows from Proposition \ref{Commensurability} and Theorem \ref{Nilpotent}.  In the case where $\phi$ is hyperbolic, the eigenvalues are real, and the conclusion follows from Theorem \ref{BoundedPacking}.
\end{proof}
\end{cor}
\section{\bf Coset Growth of $H$ in $P_{2,\phi}$}\label{CosetGr}
\subsection{\bf Growth of Groups and Coset Growth}
\subsubsection{\bf Growth of Groups}
In the course of proving bounded packing in $P_{n,\phi}$, we have developed the tools for studying the growth of a certain metric space.  Recall that given a metric measure space $(X,d,\mu)$, one defines the \textit{volume growth} to be the function $f_{x_0}(r)=\mu(B(x_0,r))$, where $B(x_0,r)$ is the open ball of radius $r$ centered at the point $x_0\in X$.  In geometric group theory, the metric space is a finitely generated group $G$ with the word metric, the measure on $G$ is the counting measure, and the basepoint is the identity element $1\in G$.  In other words, the function $f(r)$ counts the number of elements which can be expressed as a word of length at most $r$ in the generators.  In order to make the definition independent of the generating set, we introduce the following equivalence relation on the set of non-decreasing functions defined on $\mathbb{R}_{\geq 0}$: $\alpha\sim \beta$ if there exists $C>0$ such that $\beta\left(\dfrac{r}{C}\right)\leq \alpha(r)\leq \beta\left(Cr\right)$.  We now define the \textit{growth rate} of $G$ to be the equivalence class of the growth function $f(r)$.

It is important to note that the growth rate of any finitely generated group is at most exponential.  This is due to the fact that the growth rate of any quotient of the finitely generated group $G$ is bounded above by the growth rate of $G$.  On the other hand, every finitely generated group is a quotient of a finitely generated free group and these are easily shown to have exponential growth.  What else can be said of the growth rate?  The free abelian groups $\mathbb{Z}^n$ provide examples of groups of polynomial growth of any order $n$.  In \cite{Wlf}, Wolf proves that nilpotent groups have polynomial growth.  Actually, \cite{Wlf} contains a stronger result which is of interest for us:
\begin{thm}[Wolf \cite{Wlf}]\label{PolycGrthWlf}A polycyclic group is either virtually nilpotent and is thus of polynomial growth, or is not virtually nilpotent and is of exponential growth.
\end{thm}
Coupling this with a result by Milnor, see \cite{Mil}, Wolf deduces:
\begin{thm}[Wolf \cite{Wlf}]A finitely generated solvable group is either polycyclic and virtually nilpotent and is thus of polynomial growth, or has no nilpotent subgroup of finite index and is of exponential growth.
\end{thm}
Much progress has been made in the area of geometric group theory concerned with the growth of groups since the works of Wolf and Milnor.  The most notable result is, of course, Gromov's celebrated result which characterizes the finitely generated, virtually nilpotent groups, as precisely the finitely generated groups of polynomial growth:
\begin{thm}[Gromov \cite{Gr2}]A finitely generated group has polynomial growth if and only if it is virtually nilpotent.
\end{thm}
Gromov's original proof of this amazing result relies on the solution of Hilbert's Fifth Problem which, in its original version, aimed to characterize Lie groups as the topological groups which are also topological manifolds, a very powerful result in and of itself.  Since Gromov's original proof in \cite{Gr2}, simpler proofs, which do not use the solution to Hilbert's Fifth Problem, have appeared.  See, for example, the work by Kleiner in \cite{Kln}.

\subsubsection{\bf Coset Growth}\label{CosGr}
Let $G$ be a finitely generated group and $H$ be any subgroup.  Consider the space of left cosets $G/H$.  The word metric $d$ on $G$ induces a distance function on $\mathcal{P}(G)$, the set of subsets of $G$, by setting $d(X,Y)=inf\left\{d(x,y):x\in X,y\in Y\right\}$.  We define the \textit{coset growth} $f_{G,H}$ to be the growth rate of the coset space $G/H$ using the distance function $d$ instead of a true metric.  We immediately observe that $f_{G,\left\{1\right\}}$ is just the growth rate of $G$, while $f_{G,N}$ is the growth rate of the quotient group $G/N$ whenever $N$ is normal in $G$.  We shall also refer to the coset growth of $H$ in $G$ as the growth of the pair $(G,H)$.  Having an infinite coset growth of $H$ in $G$ is equivalent to $H$ not having the property of bounded packing in $G$.
Our goal in this section is to show that the growth rate of the pair $(P_{2,\phi},\left\langle t\right\rangle)$ is at most exponential, when $P_{2,\phi}$ is not nilpotent.  Theorem \ref{PolycGrthWlf} shows that the growth rate of $(P_{2,\phi},\left\{1\right\})$ is exponential, whereas the growth rate of $(P_{2,\phi},\mathbb{Z}^2)$ is polynomial of order $1$ as $P_{2,\phi}/\mathbb{Z}^2\cong\mathbb{Z}$.

In order to bound the growth rate of $(P_{2,\phi},\left\langle t\right\rangle)$ above, we need to examine the distortion of $\mathbb{Z}^2$ in $P_{2,\phi}$.  To do this, we make use of the fact that $P_{2,\phi}$ embeds as a uniform lattice of finite covolume in the Lie group $Sol$.  In fact, every polycyclic group can be virtually embedded as a lattice in a simply-connected solvable Lie group.  However, in the case of $P_{2,\phi}$ the embedding is far less abstract.
\subsection{\bf The $Sol$ Geometry}\label{Sol}
\subsubsection{\bf Definition and Basic Properties of $Sol$}
The Lie group $Sol$ is defined to be the semidirect product of $\mathbb{R}^2\rtimes\mathbb{R}$, where a generator $t\in\mathbb{R}$ acts on $\mathbb{R}^2$ via $
 \begin{pmatrix}
  e^t & 0 \\
   0 & e^{-t}
 \end{pmatrix}
$.  We identify $Sol$ as a set with $\mathbb{R}^3=\left\{(x,y,z):(x,y)\in\mathbb{R}^2, z\in\mathbb{R}\right\}$ and thereby obtain a chart on the underlying manifold.  We shall use this chart to introduce a left-invariant Riemannian metric on $Sol$ as follows: First, we metrically identify the normal copy of $\mathbb{R}^2$ with $\mathbb{E}^2$.  Next, we transport the metric $ds^2_{z=0}=dx^2+dy^2+dz^2$ on $\mathbb{R}^2\times\left\{0\right\}$ to $\mathbb{R}^2\times\left\{t\right\}$ for arbitrary $t\in\mathbb{R}$ using the multiplication map $m_t:Sol\rightarrow Sol$ given by $m_t(x,y,z)=(0,0,t)(x,y,z)$.  This means that we define $ds^2_{z=t}$ in such a way that it satisfies $ds^2_{z=t}(Dm_t(\textbf{v}))=ds^2_{z=0}(\textbf{v})$ for every $\textbf{v}\in T_p(\mathbb{R}^2)\subseteq T_p(Sol)$.  It is easy to show that the resulting left-invariant Riemannian metric is given by $ds^2=e^{2z}dx^2+e^{-2z}dy^2+dz^2$.

We shall denote the associated path metric on $Sol$ by $d_S$.  The embedding of $\mathbb{E}^2$ into $Sol$ via the map $(x,y)\mapsto (x,y,z)$ endows $\mathbb{R}^2\subseteq Sol$ with the Euclidean metric which we shall denote by $d_{\mathbb{R}^2}$.  In order to produce the promised bound for the growth rate of the pair $(P_{2,\phi},\left\langle t\right\rangle)$, we need the following result which provides a lower bound for the distance in $Sol$ between two points in $\mathbb{R}^2$ in terms of their $l^1$ distance in $\mathbb{R}^2$:
\begin{lem}\label{SolEstimate}$d_S((x_1,y_1,0),(x_2,y_2,0))\geq\max\left\{2\log|x_2-x_1|,2\log|y_2-y_1|\right\}$.
\begin{proof}See Lemma \ref{DistancesBddBelow} in \cite{NB}.
\end{proof}
\end{lem}
We now obtain a lower bound for $d_S$ in terms of the Euclidean distance:
\begin{equation}\label{estimate1}
d_{\mathbb{R}^2}((x_1,y_1),(x_2,y_2))\leq \sqrt{2}\max\left\{|x_2-x_1|,|y_2-y_1|\right\}
\end{equation}
\begin{equation}
\log d_{\mathbb{R}^2}((x_1,y_1),(x_2,y_2))\leq\max\left\{\log|x_2-x_1|,\log|y_2-y_1|\right\}+\log\frac{1}{\sqrt{2}}
\end{equation}
hence by Lemma \ref{SolEstimate}
\begin{equation}
d_S((x_1,y_1,0),(x_2,y_2,0))\geq 2\log d_{\mathbb{R}^2}((x_1,y_1,0),(x_2,y_2,0))+\log{2}
\end{equation}
This estimate bounds the distortion of $\mathbb{R}^2$ in $Sol$ and our next goal will be to quasify it to apply to the discrete case of $P_{2,\phi}$.  Before we do that however, we provide an explicit embedding of the non-nilpotent group $P_{2,\phi}$ as a uniform lattice in $Sol$ following \cite{NB}.
\subsubsection{\bf The Embedding of $P_{2,\phi}$ in $Sol$}
First, we embed $\mathbb{Z}^2$ as the standard integer lattice in $\mathbb{R}^2$.  Let the eigenvalues of $\phi$ be $0<\lambda^{-1}<1<\lambda$, and let $\textbf{v}_-$ and $\textbf{v}_+$ be eigenvectors for the corresponding eigenvalues.  Let $f:\mathbb{R}^2\rightarrow\mathbb{R}^2$ be the linear map defined by $\textbf{v}_-\mapsto e_1$ and $\textbf{v}_+\mapsto e_2$.  Now, we define a map $\psi:P_{2,\phi}\rightarrow Sol$ by $\psi(e_1^k e_2^l t^q)=(f(k,l),q\log(\lambda))$, where $e_1,e_2$ and $t$ are the standard generators of $P_{2,\phi}$ described in Section \ref{PolyGrps}, and $k,l,q\in\mathbb{Z}$.  It is not hard to show that the map $\psi$ is a monomorphism with a discrete image and that $\psi(P_{2,\phi})\backslash Sol$ is compact.
\subsubsection{\bf Bounding the Growth of $(P_{2,\phi},\left\langle t\right\rangle)$}
By uniformity of the embedding, we have the following quasi-isometry inequality:
\begin{equation}\label{DiscreteEstimate}
d_{P_{2,\phi}}((x_1,y_1,0),(x_2,y_2,0))\geq Cd_S((x_1,y_1,0),(x_2,y_2,0)).
\end{equation}
Now, passing from the continuous metric on $\mathbb{R}^2$ to the word metric on $\mathbb{Z}^2$ we have
\begin{equation}
d_{\mathbb{R}^2}((x_1,y_1,0),(x_2,y_2,0))\geq\frac{1}{\sqrt{2}}d_{\mathbb{Z}^2}((x_1,y_1,0),(x_2,y_2,0))
\end{equation}
and finally putting all these estimates together we get:
\begin{equation}\label{DistortionBound}
d_{P_{2,\phi}}((x_1,y_1,0),(x_2,y_2,0))\geq C\log d_{\mathbb{Z}^2}((x_1,y_1,0),(x_2,y_2,0))+A
\end{equation}
\begin{thm}\label{Cogrowth}The coset growth $f_H(r)$ of $H=\left\langle t\right\rangle$ in $P_{2,\phi}$ is at most exponential.
\begin{proof}Given $r>0$, Lemma \ref{DistancesBddBelow} shows that there exists $N\sim 2 \left(\pi r^2\right)^2$ such that any $N$ distinct cosets of $H$, say $a_1H,...,a_NH$, contains a pair $aH,bH$ with $d_{\mathbb{Z}^2}(\phi^{-k}(a),\phi^{-k}(b))>r$.  Now, we have $d_{P_{2,\phi}}(at^k,bt^k)=d_{P_{2,\phi}}(\phi^{-k}(a),\phi^{-k}(b))$, which combined with (\ref{DistortionBound}) gives $d_{P_{2,\phi}}(at^k,bt^k)>C\log r+A$.  Arguing as in the proof of Lemma \ref{BddPackoft}, we conclude that $d_{P_{2,\phi}}(aH,bH)\geq\frac{1}{2}(C\log r+A)$.  Now, setting $R=\frac{1}{4}(C\log r+A)-1$, we conclude that the open ball of radius $R$ in $P_{2,\phi}/H$ has at most $N\sim B\alpha^r$ elements, where $B$ and $\alpha$ depend only on $A$ and $C$.  
\end{proof}
\end{thm}











\bibliographystyle{plain}

\end{document}